\documentclass[11pt,bezier]{article}
\usepackage{amsmath,amssymb,amsfonts}
\usepackage{currfile}
\usepackage{graphicx,xcolor}
\usepackage[colorlinks,linkcolor=red,citecolor=blue]{hyperref}
\usepackage{pgf,tikz}
\usepackage{mathrsfs}
\usetikzlibrary{arrows}

\newtheorem{thm}{Theorem}[section]
\newtheorem{defi}[thm]{Definition}
\newtheorem{remark}[thm]{Remark}
\newtheorem{cor}[thm]{Corollary}
\newtheorem{obs}[thm]{Observation}
\newtheorem{prop}[thm]{Proposition}
\newtheorem{lem}[thm]{Lemma}
\newtheorem{example}[thm]{Example}

\newenvironment{proof}[1]{{\bf Proof. }{\rm #1}}{\hfill $\rule {2mm}{2mm}$\\}

\newtheorem{prelem}{{\bf Theorem}}

\newenvironment{oldtheorem}{\begin{prelem}{\hspace{-0.5
em}{\bf}}}{\end{prelem}}

\begin{document}

 \title{A study of $4-$cycle systems  }
\author{
{\sc B. Bagheri Gh.${}^{a}$\footnote{behrooz@ac.tuwien.ac.at}, M. Khosravi${}^{b}$\footnote{
  		khosravi$_-$m@uk.ac.ir (Corresponding author)}, E. S. Mahmoodian${}^{c}$\footnote{emahmood@sharif.edu} and S. Rashidi${}^{d}$\footnote{
  		saeedeh.rashidi@uk.ac.ir}}}
  	\date{\bf{ }}
\maketitle
\vspace{-1cm}
  \begin{center}
  	${}^a$
{\small \it Algorithms and Complexity Group}\\
  	{\small \it Vienna University of Technology}\\
  	{\small  \it Favoritenstrasse 9-11,}
  \vspace*{5mm}
  	{\small \it 1040 Vienna, Austria }\\
  	  	${}^b$
  	{\small \it Department of Pure Mathematics}\\
  	{\small \it Faculty of Mathematics and Computer}\\
  	{\small  \it Shahid Bahonar University of Kerman,}
  	\vspace*{5mm}
  	{\small \it Kerman, Iran} \\
  	${}^c$
{\small \it Department of Mathematical Sciences} \\
{\small \it Sharif University of Technology} \\
{\small \it P.O. Box 11155-9415,}
\vspace*{5mm}
{\small \it Tehran, Iran} \\
   ${}^d$
   {\small \it Department of Applied Mathematics}\\
  	{\small \it Faculty of Mathematics and Computer}\\
  	{\small  \it Shahid Bahonar University of Kerman,}
  	\vspace*{5mm}
  	{\small \it Kerman, Iran} \\
  \end{center}

% ---------------------------------------------------------------------------------------------------------------------------------------------------------------------------------------

 % -----------------------------------------------------------------------------------------------------------------------------------------------------------------------------------
 \begin{abstract}
A $4-$cycle system is a partition of the edges of the complete graph $K_n$ into $4-$cycles.
Let ${ C}$ be a collection of cycles of length 4 whose edges partition the edges of   $K_n$.
A set of 4-cycles $T_1 \subset C$
is called a
4-cycle trade
if there exists a set
 $T_2$ of edge-disjoint
 4-cycles on the same vertices,
 such that
$
({C} \setminus T_1)\cup T_2$
also
is
a collection of cycles of length 4 whose edges partition the edges  of
$K_n$.

 We study  $4-$cycle trades of volume two (double-diamonds) and three and show that the set of all 4-CS(9) is connected with respect of trading with  trades of volume 2 (double-diamond) and 3.

In addition, we present a full rank matrix whose null-space is containing trade-vectors.
\end{abstract}
\hspace*{-2.7mm} {\bf MSC:} {\sf 05B30; 05B05}\\
\hspace*{-2.7mm} {\bf KEYWORDS:} { \sf $4-$cycle system, double-diamond, $4-$cycle trade.}
% ---------------------------------------------------------------------------------------------------------------------------------------------------------------------------------------------------------------------------------------
\maketitle
% ----------------------------------------------------------------------------------------------------------------------------------------------------------------------------------------------------------------------------------------
%%%%%%%%%%%%%%%%%%%%%%%%%%%%%%%%%%%%%%%%%%%%%%%%%%%%%%%%%%%%%%%%%%%%%%%%%%%%%%%%%%%%%%%%%%%
%-----------------------------------------------------------------------------------------------------------------------------------------------------------------------------------------------------------------------------------------
\section{Introduction}
Let $G$ be a graph with the vertex set $V(G)$ and the edge set
$E(G)$.  An $n$-cycle $(v_1,\ldots,v_n)$ (or briefly $v_1v_2\dots
v_n$) in a graph $G$, is a subgraph of $G$ which consists of $n$
distinct vertices, $v_i$'s, and edge set
$\{\{v_1,v_2\},\{v_2,v_3\},\ldots,\{v_n,v_1\}\}$. A {\sf
$4-$cycle system} of order $n$, denoted by $4-{\rm CS}(n)$, is a
collection  of cycles of length $4$ whose edges partition the edges
of $K_n$.\\
The following well-known theorem states that for which values of $n$ a $4-{\rm CS}(n)$ exists.
\begin{oldtheorem}~{\rm(\cite[Page \ 266]{MR1392993})}
    A necessary
    and sufficient condition for the existence of a $4-{\rm CS}(n)$ is that
    $n \equiv 1 (\hspace{-2mm}\mod 8)$.
\end{oldtheorem}
An interesting problem in combinatorics is that  whether  there can be defined some ``moves"
(using ``trades'' of small volume or something else) between different elements of a
class of combinatorial objects with the same parameters, such as Latin squares, Steiner
triple systems, etc. These moves must be such a way that each element has chance to be produced by these moves.

By simulating an ergodic Markov chain whose stationary distribution is
uniform over the space of $n \times n$ Latin squares, Mark T. Jacobson and Peter
Matthews~\cite{MR1410617}, have discussed elegant method by which they generate Latin
squares with a uniform distribution (approximately). The central issue is
 finding the moves that connect the squares. Aryapoor and Mahmoodian  obtained
a short proof for that, by using Latin trades of volume four (i.e. intercalates) \cite{MR2883525}. \\
There does not exist a known move between Steiner triple systems as yet. Steiner triple systems are
$3-$cycle systems. Here, we investigate the $4-$cycle systems and $4-$cycle trades.

\begin{defi}
Let $T_1$ be a set of edge-disjoint $4-$cycles on the vertex set
$\{1,\ldots,v\}$. Then $T_1$ is called a {\sf $4-$cycle trade},
if there exists a set, $T_{2}$, of edge-disjoint $4-$cycles on the same vertex set
$\{1,\ldots,v\}$, such that $T_{1}\cap T_2=\emptyset$ and
$\bigcup_{C\in T_1}E(C)=\bigcup_{C\in T_{2}}E(C)$.
\end{defi}
 We call $T_{2}$ a {\sf disjoint mate} of $T_1$ and the pair $\mathcal{T}=(T_1,T_{2})$ is called a
{\sf $4-$cycle bitrade} of {\sf{volume}} $s=|T_1|$ and {\sf foundation} $v=|\bigcup_{C\in T_1}V(C)|$. Here for them  we use the term ``trade'' for short.

A {\sf $\mu-$way $4-$cycle trade} is a collection of $\mu$ disjoint
collections $\{T_1,\dots,T_{\mu}\}$ such that $(T_i,T_j)$ forms a
$4-$cycle bitrade for each $i\neq j$.

\section{$4-$cycle trades of small volumes}
\subsection{Trades of volume 2}
\begin{figure}[!h]\centering
\begin{tikzpicture}
\clip(-4.3,-2.48) rectangle (7.06,6.3); \draw (-3.46,5.36)--
(-2.26,5.4); \draw (-2.26,5.4)-- (-2.24,4.34); \draw
(-2.24,4.34)-- (-3.46,4.32); \draw (-3.46,4.32)-- (-3.46,5.36);
\draw (-1.38,5.34)-- (-0.32,5.34); \draw (-0.32,5.34)--
(-0.3,4.38); \draw (-0.3,4.38)-- (-1.4,4.36); \draw (-1.4,4.36)--
(-1.38,5.34); \draw (2.72,5.6)-- (3.5,4.84); \draw (3.5,4.84)--
(2.72,4.); \draw (2.72,4.)-- (1.96,4.76); \draw (1.96,4.76)--
(2.72,5.6); \draw (3.5,4.84)-- (4.3,5.72); \draw (4.3,5.72)--
(5.16,4.96); \draw (5.16,4.96)-- (4.4,4.06); \draw (4.4,4.06)--
(3.5,4.84); \draw (-3.08,2.42)-- (-3.06,1.3); \draw (-3.06,1.3)--
(-1.88,1.3); \draw (-1.88,1.3)-- (-1.88,2.48); \draw
(-1.88,2.48)-- (-3.08,2.42); \draw (-1.88,2.48)-- (-3.06,1.3);
\draw (-3.06,1.3)-- (-1.8,0.34); \draw (-1.88,2.48)--
(-0.74,1.46); \draw (-0.74,1.46)-- (-1.8,0.34); \draw
(3.58,2.42)-- (5.38,1.48); \draw (5.38,1.48)-- (3.64,0.54); \draw
(3.58,2.42)-- (1.96,1.4); \draw (1.96,1.4)-- (3.64,0.54); \draw
(3.64,0.54)-- (3.1,1.48); \draw (3.1,1.48)-- (3.58,2.42); \draw
(3.58,2.42)-- (4.1,1.5); \draw (4.1,1.5)-- (3.64,0.54);
\begin{scriptsize}
\draw [fill=black] (-3.46,5.36) circle (2.5pt); \draw
[fill=black] (-2.26,5.4) circle (2.5pt); \draw [fill=black]
(-3.46,4.32) circle (2.5pt); \draw [fill=black] (-2.24,4.34)
circle (2.5pt); \draw [fill=black] (-1.38,5.34) circle (2.5pt);
\draw [fill=black] (-0.32,5.34) circle (2.5pt); \draw
[fill=black] (-0.3,4.38) circle (2.5pt); \draw [fill=black]
(-1.4,4.36) circle (2.5pt); \draw [fill=black] (2.72,5.6) circle
(2.5pt); \draw [fill=black] (3.5,4.84) circle (2.5pt); \draw
[fill=black] (2.72,4.) circle (2.5pt); \draw [fill=black]
(1.96,4.76) circle (2.5pt); \draw [fill=black] (4.3,5.72) circle
(2.5pt); \draw [fill=black] (5.16,4.96) circle (2.5pt); \draw
[fill=black] (4.4,4.06) circle (2.5pt); \draw [fill=black]
(-3.08,2.42) circle (2.5pt); \draw [fill=black] (-3.06,1.3)
circle (2.5pt); \draw [fill=black] (-1.88,1.3) circle (2.5pt);
\draw [fill=black] (-1.88,2.48) circle (2.5pt); \draw
[fill=black] (-1.8,0.34) circle (2.5pt); \draw [fill=black]
(-0.74,1.46) circle (2.5pt); \draw [fill=black] (3.58,2.42)
circle (2.5pt); \draw [fill=black] (5.38,1.48) circle (2.5pt);
\draw [fill=black] (3.64,0.54) circle (2.5pt); \draw
[fill=black] (1.96,1.4) circle (2.5pt); \draw [fill=black]
(3.1,1.48) circle (2.5pt); \draw [fill=black] (4.1,1.5) circle
(2.5pt);
\end{scriptsize}
\end{tikzpicture}
 \vspace{-2cm}\caption{All possible configurations of two $4-$cycles.}\label{two4cycles}
\end{figure}
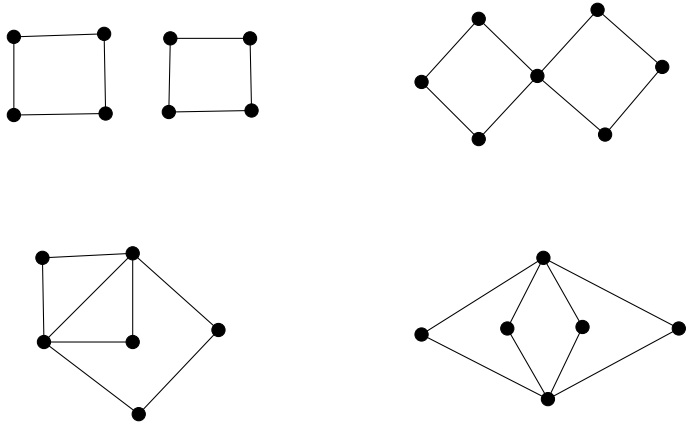

There are
just four possible configurations of two disjoit $4-$cycles (Figure~\ref{two4cycles}).
%\vspace{-2cm}
It is easily seen that only the last configuration in Figure~\ref{two4cycles} may form a
bitrade called a {\sf double-diamond} $\mathcal{D}=(D_1,D_2)$ in which  $D_1=  \{(1,2,3,4),(1,5,3,6)\}$ and $D_2= \{(1,2,3,5)$, $(1,4,3,6)\}$, where the
vertices are labeled $\{1,\ldots,6\}$.

 We have the following $3-$way $4-$cycle trade of volume $2$ and foundation $6$.
\begin{table}[h!]\begin{center}
        \begin{tabular}{c|c|c}
            $T_1$&$T_2$&$T_3$\\
            \hline
            $1234$&$1235$&$1236$\\
            $1536$& $1436$&$1435$\\
        \end{tabular}
\caption{The $3-$way $4-$cycle trades of
volume $2$.}
\label{table3w4ct2}\end{center}
\end{table}

   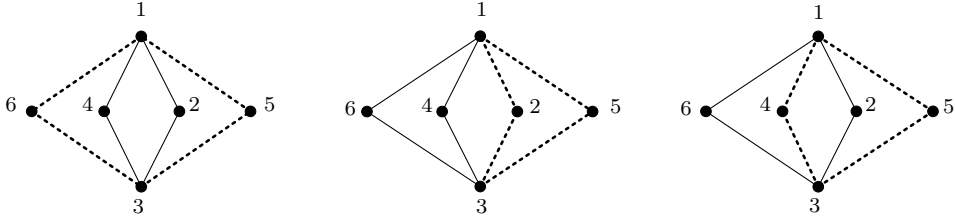
\begin{figure}[!h]\centering
\begin{tikzpicture}[line cap=round,line join=round,>=triangle 45,x=1.0cm,y=1.0cm]
\clip(-1.,-0.5) rectangle (12.,2.5); \draw (1.,2.)--
(1.506727272727275,1.0027272727272722); \draw
(1.506727272727275,1.0027272727272722)-- (1.,0.); \draw (1.,0.)--
(0.5067272727272748,1.0027272727272722); \draw
(0.5067272727272748,1.0027272727272722)-- (1.,2.); \draw
[line width=1.pt,dotted] (1.,2.)-- (-0.4569090909090887,1.0027272727272722);
\draw [line width=1.pt,dotted] (-0.4569090909090887,1.0027272727272722)--
(1.,0.); \draw [line width=1.pt,dotted] (1.,0.)--
(2.45218181818182,1.0027272727272722); \draw [line width=1.pt,dotted]
(2.45218181818182,1.0027272727272722)-- (1.,2.); \draw
(5.506727272727274,2.002727272727271)-- (5.,1.); \draw (5.,1.)--
(5.506727272727274,0.0027272727272733446); \draw [line width=1.pt,dotted]
(5.506727272727274,0.0027272727272733446)-- (6.,1.); \draw
[line width=1.pt,dotted] (6.,1.)-- (5.506727272727274,2.002727272727271);
\draw (5.506727272727274,2.002727272727271)-- (4.,1.); \draw
(4.,1.)-- (5.506727272727274,0.0027272727272733446); \draw
[line width=1.pt,dotted] (5.506727272727274,0.0027272727272733446)--
(7.,1.); \draw [line width=1.pt,dotted]
(5.506727272727274,2.002727272727271)-- (7.,1.); \draw
[line width=1.pt,dotted] (10.,2.)-- (9.524909090909093,1.0027272727272722);
\draw (10.,2.)-- (10.506727272727275,1.0027272727272722); \draw
(10.506727272727275,1.0027272727272722)-- (10.,0.); \draw
[line width=1.pt,dotted] (10.,0.)-- (9.524909090909093,1.0027272727272722);
\draw (10.,2.)-- (8.506727272727275,1.0027272727272722); \draw
(8.506727272727275,1.0027272727272722)-- (10.,0.); \draw
[line width=1.pt,dotted] (10.,0.)-- (11.524909090909093,1.0027272727272722);
\draw [line width=1.pt,dotted] (11.524909090909093,1.0027272727272722)--
(10.,2.);
\begin{scriptsize}
\draw [fill=black] (1.,2.) circle (2.0pt); \draw[color=black]
(0.9976363636363658,2.3572727272727256) node {$1$}; \draw
[fill=black] (1.506727272727275,1.0027272727272722) circle
(2.0pt); \draw[color=black]
(1.7067272727272746,1.1027272727272723) node {$2$}; \draw
[fill=black] (1.,0.) circle (2.0pt); \draw[color=black]
(0.9612727272727294,-0.26090909090909) node {$3$}; \draw
[fill=black] (0.5067272727272748,1.0027272727272722) circle
(2.0pt); \draw[color=black]
(0.2885454545454567,1.1027272727272723) node {$4$}; \draw
[fill=black] (-0.4569090909090887,1.0027272727272722) circle
(2.0pt); \draw[color=black]
(-0.7296363636363614,1.1027272727272723) node {$6$}; \draw
[fill=black] (2.45218181818182,1.0027272727272722) circle (2.0pt);
\draw[color=black] (2.706727272727275,1.1027272727272723) node
{$5$}; \draw [fill=black] (5.506727272727274,2.002727272727271)
circle (2.0pt); \draw[color=black]
(5.524909090909093,2.3572727272727256) node {$1$}; \draw
[fill=black] (5.,1.) circle (2.0pt); \draw[color=black]
(4.797636363636366,1.084545454545454) node {$4$}; \draw
[fill=black] (5.506727272727274,0.0027272727272733446) circle
(2.0pt); \draw[color=black] (5.524909090909093,-0.26090909090909)
node {$3$}; \draw [fill=black] (6.,1.) circle (2.0pt);
\draw[color=black] (6.234,1.066363636363636) node {$2$}; \draw
[fill=black] (4.,1.) circle (2.0pt); \draw[color=black]
(3.779454545454547,1.0481818181818177) node {$6$};
%\draw[color=black] (4.634,0.43) node {$t$};
\draw [fill=black] (7.,1.) circle (2.0pt); \draw[color=black]
(7.306727272727275,1.066363636363636) node {$5$}; \draw
[fill=black] (10.,2.) circle (2.0pt); \draw[color=black]
(10.006727272727275,2.320909090909089) node {$1$}; \draw
[fill=black] (9.524909090909093,1.0027272727272722) circle
(2.0pt); \draw[color=black] (9.297636363636365,1.1027272727272723)
node {$4$}; \draw [fill=black]
(10.506727272727275,1.0027272727272722) circle (2.0pt);
\draw[color=black] (10.697636363636365,1.1027272727272723) node
{$2$}; \draw [fill=black] (10.,0.) circle (2.0pt);
\draw[color=black] (9.95218181818182,-0.24272727272727188) node
{$3$}; \draw [fill=black] (8.506727272727275,1.0027272727272722)
circle (2.0pt); \draw[color=black]
(8.243090909090912,1.0481818181818177) node {$6$}; \draw
[fill=black] (11.524909090909093,1.0027272727272722) circle
(2.0pt); \draw[color=black] (11.734,1.066363636363636) node {$5$};
\end{scriptsize}
\end{tikzpicture}
       \caption{A $3-$way $4-$cycle trade of volume $2$ (double-diamond).}\label{15}
    \end{figure}

The following
theorem is proved by two different  approaches, for example see~\cite{MR2195316} and a reference in there.
\begin{oldtheorem}~{\rm(\cite{MR2195316})}\label{diamond}
For each $n=8k+1$, there exists a  $4-{\rm CS}(n)$   which does not contain
any double-diamond.
\end{oldtheorem}
So if we use double-diamond trades as our ``moves", we will not be able to produce ``all" $4-{\rm CS}(n)$s.

\subsection{Trades of volume 3}

    Let   $(T_1,T_2)$ be a
    {$4-$cycle trade} of {volume} $s$ and {foundation}
    $v$. Then $G$,  the {\sf graph of this trade} is a graph obtained by the union of the edges in all $4-$cycles
    in $T_1$ which is the same if we use the cycles of $T_2$. The graph $G$ has $v$ vertices and $4s$ edges.
    All the vertices of $G$ have non-negative even  degrees, say $d_i$, $i=1,\ldots,v$.
    We denote the number of vertices
    of degree $j$ in $G$ by $x_j$.   We have the following equations,
    \begin{equation}
        \label{sum of degrees}
        x_1+x_2+\cdots+x_{v-1}=v,  \ \ \ \
        x_1+2x_2+\cdots+(v-1)x_{v-1}=8s.
    \end{equation}
    Note that  $x_{2k+1}=0$, for each non-negative number $k$.
\begin{obs} \label{obs}
 If a bridgeless graph $G$ has two adjacent vertices of degree $2$, then $G$ has no two disjoint $4-$cycle decompositions.
\end{obs}
By Observation~\ref{obs}, for every  $4-$cycle trade of volume $s$ we have $x_2 \le 2s$.
\begin{thm}
There is (up to isomorphism) just one $4-$cycle bitrade of volume $3$   and
foundation $6$. This bitrade can be extended to two $3-$way  $4-$cycle trades.
\end{thm}
\begin{proof}
      Equations~(\ref{sum of degrees}) for the graph of the desired trade will be as follows,
    \begin{equation}
        \label{s=3 f=6}
        x_2+x_4=6, \ \  \ \
        2x_2+4x_{4}=24.
    \end{equation}
    The only solution for the degrees of $G$ is: $ x_2=0$  and $x_4=6$.
    The degree sequence of $G$ will be: $<4,4,4,4,4,4>$. Therefore, $\overline{G}$ is a $1-$factor, say
    $F=\{\{1,2\},\{3,4\},\{5,6\}     \}$, the non-edges of $G$ and
    $G=K_6 \setminus F$.
    Note that every edge decomposition of $G$ into $4-$cycles has at least one cycle $wxyz$ such that $\{\{w,y\},\{x,z\}\}\subset F$.
  To see this, consider one of the cycles with vertices 1 and 2. If it doesn't meet this condition, then it should be for example 1325. Now, the cycle that using other edges with endpoint 3 should be 5364, which has the desired property. Using this property, every trade is isomorphic to $(T_1,T_2)$, which are introduced in Table \ref{table3w4ct3}. The third decomposition can be of the form $T_3$ or $T_4$. It can easily seen that $(T_1,T_2,T_3)$ and $(T_1,T_2,T_4)$ are not isomorphic.
 See Figure~\ref{Fig:s=3 f=6}.
  \begin{center}
  \end{center}
\begin{figure}[!h]\centering
\begin{tikzpicture}[line cap=round,line join=round,>=triangle 45,x=1.0cm,y=1.0cm]
\clip(-0.5,-0.5) rectangle (15.5,2.5);
\draw [line width=1.pt,dash pattern=on 4pt off 4pt] (1.14,2.)-- (0.24,0.54);
\draw [line width=1.pt,dash pattern=on 4pt off 4pt] (1.14,2.)-- (0.52,0.);
\draw [line width=1.pt,dotted] (1.14,2.)-- (2.68,0.62);
\draw [line width=1.pt,dotted] (1.14,2.)-- (2.42,0.);
\draw (0.24,0.54)-- (1.86,2.);
\draw [line width=1.pt,dash pattern=on 4pt off 4pt] (0.24,0.54)-- (2.68,0.62);
\draw (0.24,0.54)-- (2.42,0.);
\draw (0.52,0.)-- (1.86,2.);
\draw [line width=1.pt,dash pattern=on 4pt off 4pt] (0.52,0.)-- (2.68,0.62);
\draw (0.52,0.)-- (2.42,0.);
\draw [line width=1.pt,dotted] (1.86,2.)-- (2.42,0.);
\draw [line width=1.pt,dotted] (1.86,2.)-- (2.68,0.62);
\draw [line width=1.pt,dotted] (5.,2.)-- (4.1,0.52);
\draw [line width=1.pt,dash pattern=on 4pt off 4pt] (5.,2.)-- (4.54,0.);
\draw [line width=1.pt,dash pattern=on 4pt off 4pt] (5.,2.)-- (6.42,0.);
\draw [line width=1.pt,dotted] (5.,2.)-- (6.8,0.58);
\draw [line width=1.pt,dotted] (4.1,0.52)-- (5.68,2.);
\draw (4.1,0.52)-- (6.8,0.58);
\draw (4.1,0.52)-- (6.42,0.);
\draw [line width=1.pt,dash pattern=on 4pt off 4pt] (4.54,0.)-- (5.68,2.);
\draw (4.54,0.)-- (6.8,0.58);
\draw (4.54,0.)-- (6.42,0.);
\draw [line width=1.pt,dash pattern=on 4pt off 4pt] (5.68,2.)-- (6.42,0.);
\draw [line width=1.pt,dotted] (5.68,2.)-- (6.8,0.58);
\draw (9.,2.)-- (8.16,0.48);
\draw (9.,2.)-- (8.46,-0.02);
\draw [line width=1.pt,dotted] (9.,2.)-- (10.42,0.02);
\draw [line width=1.pt,dotted] (9.,2.)-- (10.76,0.58);
\draw (8.16,0.48)-- (9.8,2.);
\draw [line width=1.pt,dash pattern=on 4pt off 4pt] (8.16,0.48)-- (10.76,0.58);
\draw [line width=1.pt,dash pattern=on 4pt off 4pt] (8.16,0.48)-- (10.42,0.02);
\draw (8.46,-0.02)-- (9.8,2.);
\draw [line width=1.pt,dotted] (8.46,-0.02)-- (10.42,0.02);
\draw [line width=1.pt,dotted] (8.46,-0.02)-- (10.76,0.58);
\draw [line width=1.pt,dash pattern=on 4pt off 4pt] (9.8,2.)-- (10.76,0.58);
\draw [line width=1.pt,dash pattern=on 4pt off 4pt] (9.8,2.)-- (10.42,0.02);
\draw (13.,2.)-- (12.16,0.48);
\draw (13.,2.)-- (12.46,-0.02);
\draw [line width=1.pt,dotted] (13.,2.)-- (14.42,0.02);
\draw [line width=1.pt,dotted] (13.,2.)-- (14.76,0.58);
\draw (12.16,0.48)-- (13.8,2.);
\draw [line width=1.pt,dash pattern=on 4pt off 4pt]  (12.46,-0.02)-- (14.76,0.58);
\draw [line width=1.pt,dash pattern=on 4pt off 4pt]  (12.46,-0.02)-- (14.42,0.02);
\draw (12.46,-0.02)-- (13.8,2.);
\draw [line width=1.pt,dotted]  (12.16,0.48)-- (14.42,0.02);
\draw [line width=1.pt,dotted]  (12.16,0.48)-- (14.76,0.58);
\draw [line width=1.pt,dash pattern=on 4pt off 4pt] (13.8,2.)-- (14.76,0.58);
\draw [line width=1.pt,dash pattern=on 4pt off 4pt] (13.8,2.)-- (14.42,0.02);
\begin{scriptsize}
\draw [fill=black] (1.14,2.) circle (2.0pt);
\draw[color=black] (1.04,2.3) node {$1$};
\draw [fill=black] (0.24,0.54) circle (2.0pt);
\draw[color=black] (-0.06,0.61) node {$6$};
\draw [fill=black] (0.52,0.) circle (2.0pt);
\draw[color=black] (0.54,-0.29) node {$5$};
\draw [fill=black] (2.68,0.62) circle (2.0pt);
\draw[color=black] (2.92,0.69) node {$3$};
\draw [fill=black] (2.42,0.) circle (2.0pt);
\draw[color=black] (2.56,-0.17) node {$4$};
\draw [fill=black] (1.86,2.) circle (2.0pt);
\draw[color=black] (2.,2.3) node {$2$};
\draw [fill=black] (5.,2.) circle (2.0pt);
\draw[color=black] (4.88,2.3) node {$1$};
\draw [fill=black] (4.1,0.52) circle (2.0pt);
\draw[color=black] (3.84,0.57) node {$6$};
\draw [fill=black] (4.54,0.) circle (2.0pt);
\draw[color=black] (4.6,-0.25) node {$5$};
\draw [fill=black] (6.42,0.) circle (2.0pt);
\draw[color=black] (6.56,-0.23) node {$4$};
\draw [fill=black] (6.8,0.58) circle (2.0pt);
\draw[color=black] (7.08,0.65) node {$3$};
\draw [fill=black] (5.68,2.) circle (2.0pt);
\draw[color=black] (5.72,2.3) node {$2$};
\draw [fill=black] (9.,2.) circle (2.0pt);
\draw[color=black] (8.96,2.29) node {$1$};
\draw [fill=black] (8.16,0.48) circle (2.0pt);
\draw[color=black] (7.9,0.51) node {$6$};
\draw [fill=black] (8.46,-0.02) circle (2.0pt);
\draw[color=black] (8.46,-0.25) node {$5$};
\draw [fill=black] (10.42,0.02) circle (2.0pt);
\draw[color=black] (10.4,-0.25) node {$4$};
\draw [fill=black] (10.76,0.58) circle (2.0pt);
\draw[color=black] (10.98,0.65) node {$3$};
\draw [fill=black] (9.8,2.) circle (2.0pt);
\draw[color=black] (9.76,2.27) node {$2$};
\draw [fill=black] (13.,2.) circle (2.0pt);
\draw[color=black] (12.96,2.29) node {$1$};
\draw [fill=black] (12.16,0.48) circle (2.0pt);
\draw[color=black] (11.9,0.51) node {$6$};
\draw [fill=black] (12.46,-0.02) circle (2.0pt);
\draw[color=black] (12.46,-0.25) node {$5$};
\draw [fill=black] (14.42,0.02) circle (2.0pt);
\draw[color=black] (14.4,-0.25) node {$4$};
\draw [fill=black] (14.76,0.58) circle (2.0pt);
\draw[color=black] (14.98,0.65) node {$3$};
\draw [fill=black] (13.8,2.) circle (2.0pt);
\draw[color=black] (13.76,2.27) node {$2$};
\end{scriptsize}
\end{tikzpicture}
\caption{A $3-$way $4-$cycle trade of volume $3$ and foundation
$6$.}\label{Fig:s=3 f=6}
\end{figure}
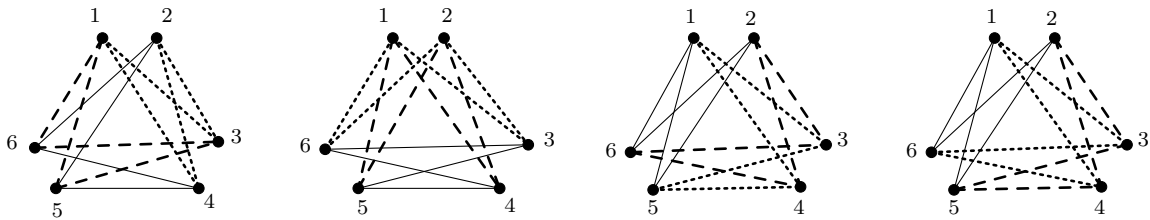    %
%\vspace{-4cm}

\begin{table}[h!]\begin{center}
        \begin{tabular}{c|c|c|c}
            $T_1$&$T_2$&$T_3$&$T_4$\\
            \hline
            $1324$&$1326$&$1354$&$1364$\\
            $1536$& $1425$&$2364$&$2354$\\
            $2546$&$3546$&$1526$&$1526$
        \end{tabular}
\caption{The $3-$way $4-$cycle trades of
volume $3$.}
\label{table3w4ct3}\end{center}
\end{table}
\vspace{-1cm}
\end{proof}

\begin{remark}
Start from
$T_1$ and replace  the two $4-$cycles $(1536)$ and $(2546)$ with $(1526)$ and
$(3546)$. Then change the two cycles $(1526)$ and $(1324)$ to $(1326)$ and
$(1425)$ which is $T_2$. Thus,  if we start from one of the $4-$cycle trade of volume $3$ and foundation $6$, by twice trade off of
double-diamonds, we can get the other one. In other
word,
 every $4-$cycle trade of volume $3$ and foundation $6$, can be obtained
from another $4-$cycle trade of volume $3$ and foundation $6$
by  trade off a  $4-$cycle bitrades of volume $2$, i.e. a double-diamond. 
\end{remark}
Since every $s\geq2$ can be written as the form $s=2n+3m$ for some $n,m\geq0$, we have the following result.
\begin{cor}
There exists $3-$way $($and therefore $2-$way$)$ $4-$cycle trade, of
volume $s$ for all $s\geq 2$ and appropriate foundation.
\end{cor}
%\begin{proof}
%The linked file is $R2$.
%\end{proof}
 \begin{thm}\label{Thm:Tv3f7}
There are  (up to isomorphism) just three  $4-$cycle bitrade of volume $3$ and foundation $7$.
\end{thm}
\begin{proof}
%
%(The linked files are $kh1$ and $B2$.)\\
    Assume  $(T_1,T_2)$ is a
    {$4-$cycle bitrade} of {volume} $3$ and {foundation}
    $7$. Then  Equations~(\ref{sum of degrees})
    for $G$, the graph of this bitrade, will be as follows,
    \begin{equation}
        \label{s=3 f=7}
        x_2+x_4+x_6=7, \ \  \ \
        2x_2+4x_{4}+6x_6=24.
    \end{equation}
    There are $3$ sets of solutions for the degrees of $G$:
    \begin{enumerate}
        \item
        $ x_2=4$,  and $x_4=1$,   $x_6=2$,
        \item
        $ x_2=3$,  and $x_4=3$,   $x_6=1$,
        \item
        $ x_2=2$,  and $x_4=5$,   $x_6=0$.
    \end{enumerate}
    The degree sequence $<d_1,d_2,d_3,d_4,d_5,d_6,d_7>$ for
    $G$ may be either one of the following,
    \begin{enumerate}
        \item
        $<2,2,2,2,4,6,6>$,
        \item
        $<2,2,2,4,4,4,6>$,
        \item
        $<2,2,4,4,4,4,4>$.
    \end{enumerate}
    We discuss each case.
    \begin{enumerate}
        \item
        $<2,2,2,2,4,6,6>$.
        This sequence is not graphical by Havel-Hakimi
        Theorem (see for example ~\cite{MR1367739}).
        \item
        {\bf\large $<2,2,2,4,4,4,6>$.} \ Let $G$ be a  realization of degree sequence $<2,2,2,4,4,4,6>$.
         Consider $V_i=\{v\in V(G)\ |\  \deg_G(v)=2i\}$, $i=1,2$. By Observation~\ref{obs}, $G(V_1)=\overline{K_3}$.
         Since  one vertex is adjacent to all other vertices, we have  $G(V_2)=K_3$. Thus, there is a unique realization of degree sequence $<2,2,2,4,4,4,6>$, see Figure~\ref{2224446}. This graph may be
        decomposed  into a $4-$cycles bitrade $\mathcal{T}^{'}=(T^{'}_1,T^{'}_2)$.
\begin{figure}[!h]\centering
 \begin{tikzpicture}[line cap=round,line join=round,>=triangle 45,x=1.0cm,y=1.0cm]
\clip(0.,-0.8) rectangle (9.,3.);
\draw [line width=1.pt,dotted] (2.,2.46)-- (3.44,1.28);
\draw (2.,2.46)-- (3.,0.5);
\draw [line width=1.pt,dash pattern=on 4pt off 4pt] (2.,2.46)-- (2.8,-0.32);
\draw [line width=1.pt,dotted] (2.,2.46)-- (2.,0.);
\draw [line width=1.pt,dash pattern=on 4pt off 4pt] (2.,2.46)-- (1.02,0.5);
\draw (2.,2.46)-- (0.72,1.28);
\draw [line width=1.pt,dotted] (3.44,1.28)-- (3.,0.5);
\draw (0.72,1.28)-- (1.02,0.5);
\draw [line width=1.pt,dash pattern=on 4pt off 4pt] (2.8,-0.32)-- (2.,0.);
\draw [line width=1.pt,dotted] (3.,0.5)-- (2.,0.);
\draw [line width=1.pt,dash pattern=on 4pt off 4pt] (1.02,0.5)-- (2.,0.);
\draw (1.02,0.5)-- (3.,0.5);
\draw [line width=1.pt,dotted] (6.92,2.38)-- (8.24,1.28);
\draw [line width=1.pt,dash pattern=on 4pt off 4pt] (6.92,2.38)-- (8.,0.52);
\draw [line width=1.pt,dash pattern=on 4pt off 4pt] (6.92,2.38)-- (7.68,-0.38);
\draw (6.92,2.38)-- (7.,0.);
\draw [line width=1.pt,dotted] (6.92,2.38)-- (5.96,0.52);
\draw (6.92,2.38)-- (5.48,1.16);
\draw [line width=1.pt,dash pattern=on 4pt off 4pt] (8.,0.52)-- (7.,0.);
\draw (7.,0.)-- (5.96,0.52);
\draw [line width=1.pt,dotted] (5.96,0.52)-- (8.,0.52);
\draw [line width=1.pt,dotted] (8.24,1.28)-- (8.,0.52);
\draw [line width=1.pt,dash pattern=on 4pt off 4pt] (7.68,-0.38)-- (7.,0.);
\draw (5.96,0.52)-- (5.48,1.16);
\begin{scriptsize}
\draw [fill=black] (2.,2.46) circle (2.0pt);
\draw[color=black] (1.94,2.73) node {$7$};
\draw [fill=black] (3.44,1.28) circle (2.0pt);
\draw[color=black] (3.66,1.35) node {$1$};
\draw [fill=black] (3.,0.5) circle (2.0pt);
\draw[color=black] (3.26,0.55) node {$4$};
\draw [fill=black] (2.8,-0.32) circle (2.0pt);
\draw[color=black] (2.88,-0.57) node {$2$};
\draw [fill=black] (2.,0.) circle (2.0pt);
\draw[color=black] (1.8,-0.21) node {$5$};
\draw [fill=black] (1.02,0.5) circle (2.0pt);
\draw[color=black] (0.76,0.43) node {$6$};
\draw [fill=black] (0.72,1.28) circle (2.0pt);
\draw[color=black] (0.4,1.41) node {$3$};
\draw [fill=black] (6.92,2.38) circle (2.0pt);
\draw[color=black] (6.9,2.7) node {$7$};
\draw [fill=black] (8.24,1.28) circle (2.0pt);
\draw[color=black] (8.5,1.43) node {$1$};
\draw [fill=black] (8.,0.52) circle (2.0pt);
\draw[color=black] (8.3,0.65) node {$4$};
\draw [fill=black] (7.68,-0.38) circle (2.0pt);
\draw[color=black] (7.96,-0.43) node {$2$};
\draw [fill=black] (7.,0.) circle (2.0pt);
\draw[color=black] (6.92,-0.27) node {$5$};
\draw [fill=black] (5.96,0.52) circle (2.0pt);
\draw[color=black] (5.68,0.49) node {$6$};
\draw [fill=black] (5.48,1.16) circle (2.0pt);
\draw[color=black] (5.22,1.21) node {$3$};
\end{scriptsize}
\end{tikzpicture}
       \caption{$\mathcal{T}^{'}=(T^{'}_1,T^{'}_2)$.}\label{2224446}
       \end{figure}
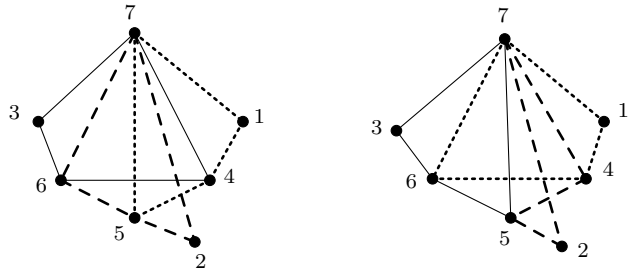
\begin{table}[h!]\begin{center}
        \begin{tabular}{c|c}
            $T^{'}_1$&$T^{'}_2$\\
            \hline
            $7145$&$7146$\\
            $7256$& $7254$\\
            $7364$&$7365$
        \end{tabular}
\caption{$\mathcal{T}^{'}=(T^{'}_1,T^{'}_2)$ a $4-$cycle bitrade of
volume $3$ and foundation $7$.}
\label{table2w4ct3}\end{center}
\end{table}

        \item
        {\bf\large$<2,2,4,4,4,4,4>$.}
         Let $G$ be a  realization of degree sequence $<2,2,4,4,4,4,4>$.
         Consider $V_i=\{v\in V(G)\ |\  \deg_G(v)=2i\}$, $i=1,2$. Let $V_1=\{v_1,v_2\}$.
By Observation~\ref{obs},  $G(V_1)=\overline{K_2}$.  Then the degree sequence of $H=G(V_2)$
            must be one of the following three  cases:
            \begin{enumerate}
                \item
                {$<2,2,4,4,4>$.} This sequence is not graphical.
                \item
                {$<2,3,3,4,4>$.} A graph with this degree sequence can be obtained from any graph with degree sequence
                $<1,2,2,3>$, which is unique. Finally $G$ is determined uniquely as in Figure~\ref{Fig:2244444''} and its bitrade
                $\mathcal{T}^{''}=(T_1^{''},T_2^{''})$  for it is given in the following.
\begin{figure}[!h]\centering
\begin{tikzpicture}[line cap=round,line join=round,>=triangle 45,x=1.0cm,y=1.0cm]
\clip(-0.5,-0.5) rectangle (6.5,2.5);
\draw [line width=1.pt,dotted] (1.,0.)-- (2.,0.);
\draw [line width=1.pt,dotted] (2.,0.)-- (2.,1.);
\draw [line width=1.pt,dotted] (2.,1.)-- (1.,1.);
\draw [line width=1.pt,dotted] (1.,1.)-- (1.,0.);
\draw [line width=1.pt,dash pattern=on 4pt off 4pt] (1.,1.)-- (0.,1.);
\draw [line width=1.pt,dash pattern=on 4pt off 4pt] (0.,1.)-- (0.,2.);
\draw [line width=1.pt,dash pattern=on 4pt off 4pt] (0.,2.)-- (1.,2.);
\draw [line width=1.pt,dash pattern=on 4pt off 4pt] (1.,2.)-- (1.,1.);
\draw (1.,2.)-- (2.,1.);
\draw (1.,2.)-- (2.,0.);
\draw (2.,1.)-- (0.,2.);
\draw [shift={(3.,3.)}] plot[domain=3.4633432079864352:4.3906384259880475,variable=\t]({1.*3.1622776601683795*cos(\t r)+0.*3.1622776601683795*sin(\t r)},{0.*3.1622776601683795*cos(\t r)+1.*3.1622776601683795*sin(\t r)});
\draw [line width=1.pt,dash pattern=on 4pt off 4pt] (4.,1.)-- (5.,1.);
\draw [line width=1.pt,dotted] (5.,1.)-- (5.,2.);
\draw (5.,2.)-- (4.,2.);
\draw [line width=1.pt,dash pattern=on 4pt off 4pt] (4.,2.)-- (4.,1.);
\draw [line width=1.pt,dotted] (5.,1.)-- (5.,0.);
\draw [line width=1.pt,dotted] (5.,0.)-- (6.,0.);
\draw (6.,0.)-- (6.,1.);
\draw [line width=1.pt,dash pattern=on 4pt off 4pt] (6.,1.)-- (5.,1.);
\draw [line width=1.pt,dotted] (5.,2.)-- (6.,0.);
\draw [line width=1.pt,dash pattern=on 4pt off 4pt] (6.,1.)-- (4.,2.);
\draw (5.,2.)-- (6.,1.);
\draw [shift={(7.,3.)}] plot[domain=3.4633432079864352:4.3906384259880475,variable=\t]({1.*3.1622776601683795*cos(\t r)+0.*3.1622776601683795*sin(\t r)},{0.*3.1622776601683795*cos(\t r)+1.*3.1622776601683795*sin(\t r)});
\begin{scriptsize}
\draw [fill=black] (1.,0.) circle (2.0pt);
\draw[color=black] (1.,-0.21) node {$1$};
\draw [fill=black] (2.,0.) circle (2.0pt);
\draw[color=black] (1.98,-0.23) node {$4$};
\draw [fill=black] (2.,1.) circle (2.0pt);
\draw[color=black] (2.24,1.13) node {$6$};
\draw [fill=black] (1.,1.) circle (2.0pt);
\draw[color=black] (0.78,1.31) node {$3$};
\draw [fill=black] (0.,1.) circle (2.0pt);
\draw[color=black] (-0.28,1.09) node {$2$};
\draw [fill=black] (0.,2.) circle (2.0pt);
\draw[color=black] (-0.1,2.27) node {$5$};
\draw [fill=black] (1.,2.) circle (2.0pt);
\draw[color=black] (0.98,2.27) node {$7$};
\draw [fill=black] (4.,1.) circle (2.0pt);
\draw[color=black] (3.76,1.09) node {$2$};
\draw [fill=black] (5.,1.) circle (2.0pt);
\draw[color=black] (4.82,1.29) node {$3$};
\draw [fill=black] (5.,2.) circle (2.0pt);
\draw[color=black] (5.02,2.27) node {$7$};
\draw [fill=black] (4.,2.) circle (2.0pt);
\draw[color=black] (4.,2.27) node {$5$};
\draw [fill=black] (5.,0.) circle (2.0pt);
\draw[color=black] (5.06,-0.27) node {$1$};
\draw [fill=black] (6.,0.) circle (2.0pt);
\draw[color=black] (6.04,-0.21) node {$4$};
\draw [fill=black] (6.,1.) circle (2.0pt);
\draw[color=black] (6.27,1.1) node {$6$};
\end{scriptsize}
\end{tikzpicture}
                     %\begin{figure}[!h]
                    %\centering
                    %\includegraphics[width=10cm]{13.png}
                    \caption{$\mathcal{T}^{''}=(T_1^{''},T_2^{''})$.}\label{Fig:2244444''}
                \end{figure}
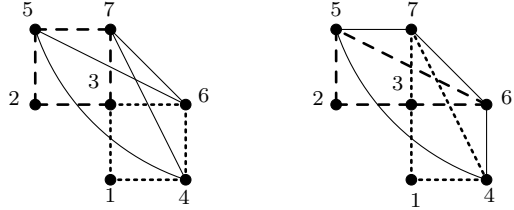
               \begin{table}[h!] \begin{center}
                \begin{tabular}{c|c}
                    $T^{''}_1$&$T^{''}_2$\\
                    \hline
                    $1364$&$1374$\\
                    $2375$& $2365$\\
                    $4567$&$4576$
                \end{tabular}
                \caption{$\mathcal{T}^{''}=(T_1^{''},T_2^{''})$ a  $4-$cycle bitrade of
volume $3$ and foundation $7$.}
\label{table2w4ct3''}\end{center}
\end{table}
                \item
                {$<3,3,3,3,4>$.} Any graph with this degree sequence can be obtained from a graph with the
                following degree sequence: $<2,2,2,2>$. And the last sequence is the degree sequence of
                a unique graph, i.e. $C_4$. So, we have two representations for
                $G$.
                \begin{enumerate}
                    \item  If $v_1$ is connected to two adjacent
                vertices of $C_4$,  then $v_1$ belongs to at most
                one cycle. Thus there is no bitrade.
                    \item  If $v_1$ is connected to
                    non-adjacent vertices of $C_4$, then $G$ can
                    be determined  uniquely as in Figure~\ref{Fig:2244444'} and a bitrade
                $\mathcal{T}^{*}=(T_1^{*},T_2^{*})$  for it is given in the following: 
\begin{figure}[!h]\centering
\begin{tikzpicture}[line cap=round,line join=round,>=triangle 45,x=1.0cm,y=1.0cm]
\clip(-0.5,-0.5) rectangle (6.5,2.5);
\draw [line width=1.pt,dotted] (1.,0.)-- (2.,0.);
\draw [line width=1.pt,dotted] (2.,0.)-- (2.,1.);
\draw [line width=1.pt,dotted] (2.,1.)-- (1.,1.);
\draw [line width=1.pt,dotted] (1.,1.)-- (1.,0.);
\draw [line width=1.pt,dash pattern=on 4pt off 4pt] (1.,1.)-- (0.,1.);
\draw [line width=1.pt,dash pattern=on 4pt off 4pt] (0.,1.)-- (0.,2.);
\draw [line width=1.pt,dash pattern=on 4pt off 4pt] (0.,2.)-- (1.,2.);
\draw [line width=1.pt,dash pattern=on 4pt off 4pt] (1.,2.)-- (1.,1.);
\draw (1.,2.)-- (2.,1.);
\draw (2.,1.)-- (0.,2.);
\draw [line width=1.pt,dash pattern=on 4pt off 4pt] (4.,1.)-- (5.,1.);
\draw (5.,1.)-- (5.,2.);
\draw (5.,2.)-- (4.,2.);
\draw [line width=1.pt,dash pattern=on 4pt off 4pt] (4.,2.)-- (4.,1.);
\draw (5.,1.)-- (5.,0.);
\draw [line width=1.pt,dotted] (5.,0.)-- (6.,0.);
\draw [line width=1.pt,dotted] (6.,0.)-- (6.,1.);
\draw [line width=1.pt,dash pattern=on 4pt off 4pt] (6.,1.)-- (5.,1.);
\draw [line width=1.pt,dash pattern=on 4pt off 4pt] (6.,1.)-- (4.,2.);
\draw [line width=1.pt,dotted] (5.,2.)-- (6.,1.);
\draw (0.,2.)-- (1.,0.);
\draw (4.,2.)-- (5.,0.);
\draw [shift={(-1.,1.)}] plot[domain=-0.46364760900080615:0.4636476090008061,variable=\t]({1.*2.23606797749979*cos(\t r)+0.*2.23606797749979*sin(\t r)},{0.*2.23606797749979*cos(\t r)+1.*2.23606797749979*sin(\t r)});
\draw [shift={(3.,1.)},line width=1.pt,dotted]  plot[domain=-0.46364760900080615:0.4636476090008061,variable=\t]({1.*2.23606797749979*cos(\t r)+0.*2.23606797749979*sin(\t r)},{0.*2.23606797749979*cos(\t r)+1.*2.23606797749979*sin(\t r)});
\begin{scriptsize}
\draw [fill=black] (1.,0.) circle (2.0pt);
\draw[color=black] (1.0036363636363639,-0.1890909090909066) node {$4$};
\draw [fill=black] (2.,0.) circle (2.0pt);
\draw[color=black] (1.9854545454545456,-0.20727272727272478) node {$2$};
\draw [fill=black] (2.,1.) circle (2.0pt);
\draw[color=black] (2.221818181818182,1.12) node {$6$};
\draw [fill=black] (1.,1.) circle (2.0pt);
\draw[color=black] (0.8036363636363639,1.2836363636363661) node {$3$};
\draw [fill=black] (0.,1.) circle (2.0pt);
\draw[color=black] (-0.2509090909090906,1.0836363636363662) node {$1$};
\draw [fill=black] (0.,2.) circle (2.0pt);
\draw[color=black] (-0.08727272727272697,2.3381818181818206) node {$5$};
\draw [fill=black] (1.,2.) circle (2.0pt);
\draw[color=black] (0.9854545454545457,2.3381818181818206) node {$7$};
\draw [fill=black] (4.,1.) circle (2.0pt);
\draw[color=black] (3.785454545454545,1.0836363636363662) node {$1$};
\draw [fill=black] (5.,1.) circle (2.0pt);
\draw[color=black] (4.84,1.2654545454545478) node {$3$};
\draw [fill=black] (5.,2.) circle (2.0pt);
\draw[color=black] (5.021818181818182,2.3381818181818206) node {$7$};
\draw [fill=black] (4.,2.) circle (2.0pt);
\draw[color=black] (4.003636363636364,2.283636363636366) node {$5$};
\draw [fill=black] (5.,0.) circle (2.0pt);
\draw[color=black] (5.0581818181818194,-0.24363636363636115) node {$4$};
\draw [fill=black] (6.,0.) circle (2.0pt);
\draw[color=black] (6.04,-0.1890909090909066) node {$2$};
\draw [fill=black] (6.,1.) circle (2.0pt);
\draw[color=black] (6.294545454545455,1.12) node {$6$};
\end{scriptsize}
\end{tikzpicture}
                %\begin{figure}[!h]
                    %\centering
                    %\includegraphics[width=10cm]{12.png}
                   \caption{$\mathcal{T}^{*}=(T^{*}_1,T^{*}_2)$.}\label{Fig:2244444'}
               \end{figure}
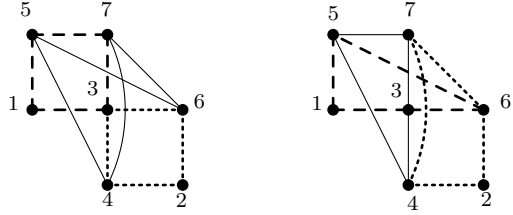

\begin{table}[h!]    \begin{center}
                \begin{tabular}{c|c}
                    $T^{*}_1$&$T^{*}_2$\\
                    \hline
                    $1375$&$1365$\\
                    $2436$& $2476$\\
                    $5476$&$4573$
                \end{tabular}
              \caption{$\mathcal{T}^{*}=(T^{*}_1,T^{*}_2)$ is a $4-$cycle bitrade of
volume $3$ and foundation $7$.}
\label{table2w4ct3*}\end{center}
\end{table}
Ofcourse, one may consider $T^*_3$ with cycles $\{5134,7365,6247\}$. But changing the places 4 and 6, $(T_1^*,T_3^*)$ will be changed to $(T_1^*,T_2^*)$.
                \end{enumerate}
            \end{enumerate}
        \end{enumerate}
        %
    %
%\vspace{-2cm}
\end{proof}
%\hspace{9cm}
%\newpage
\begin{remark}
The trade $T^{*}_1$  in the proof of Theorem~$\ref{Thm:Tv3f7}$
can be obtained from $T^{*}_2$ by twice trade off of double-diamonds.
\end{remark}
\begin{thm}
There exists precisely  one  $4-$cycle bitrades of
volume $3$ and foundation $8$. This bitrade can be extended to a $5-$way  $4-$cycle trades.
\end{thm}
\begin{proof}
    Assume $(T_1,\ldots,T_5)$ is a $5-$way
    {$4-$cycle trade} of {volume} $3$ and {foundation}
    $8$. Then  Equations~(\ref{sum of degrees})
    for $G$, the graph of this trade, will be as follows,
    \begin{equation}
        \label{5s=3 f=8}
        x_2+x_4+x_6=8, \ \  \ \
        2x_2+4x_{4}+6x_6=24.
    \end{equation}
It can be seen easily that $x_2\ge 1$.   Let $u$ be a vertex of degree $2$ in $G$ and $N_G(u)=\{x,y\}$. Since $(T_1,\ldots,T_5)$ is a $5-$way
    {$4-$cycle trade}, $x$ and $y$ have at least $5$ common neighbors except $u$ in $G$. Thus, $G$ has at least two vertices, $x$ and $y$, of degree $6$.
    So, there is exactly one set of solutions for Equation~(\ref{5s=3 f=8}) when $x_6\ge 2$ as follows
    $$x_2=6,\ x_4=0,\   x_6=2.$$
%    \begin{enumerate}
%        \item
%        $ x_2=6$,  and $x_4=0$,   $x_6=2$;
%        \item
%        $ x_2=5$,  and $x_4=2$,   $x_6=1$;
%        \item
%        $ x_2=4$,  and $x_4=4$,   $x_6=0$.
%    \end{enumerate}
    %
    Thus the degree sequence of $G$ is $<2,2,2,2,2,2,6,6>$.

    Let $V_1=\{v\in V(G)\ |\  \deg(v)=2\}$ and $V_2=V(G)\setminus V_1$. By Observation \ref{obs}, $V_1$ is an independent set; thus,
          each vertex of $V_1$ is adjacent to  all vertices of $V_2$. Since $|E(G)|=12$, the number of edges in $G[V_2]$ is
        $|E(G[V_2])|=12-|V_1||V_2|=0$. Therefore,  $G=K_{2,6}$, which can be
          decomposed
 into a $5-$way $4-$cycle trades
        (see Table~\ref{table5w4ct3}).
\begin{table}[h!]\begin{center}
\begin{tabular}{c|c|c|c|c}
$T_1$&$T_2$&$T_3$&$T_4$&$T_5$\\
\hline
$1728$&$1738$&$1748$&$1758$&$1768$\\
$3748$& $2758$&$2768$&$2748$&$2738$\\
$5768$&$4768$&$3758$&$3768$&$4758$
\end{tabular}
\caption{The $5-$way $4-$cycle trades of
volume $3$.}
\label{table5w4ct3}\end{center}
\end{table}
\end{proof}
\begin{remark}
Every trade of volume $3$ and foundation $8$ can be obtained by trades
of volume $2$ i.e. double-diamonds.
\end{remark}

\begin{thm}
There are no $4-$cycle  trade of volume $3$ and foundation more than
$8$.
\end{thm}
\begin{proof}
Suppose that $T$ is a $4-$cycle  trade of volume $3$ and foundation $v$. By  Equations~(\ref{sum of degrees}) we have
$$x_2+x_4+\cdots+x_{2k}=v,\qquad 2x_2+4x_4+\cdots+2kx_{2k}=24,$$
where $k=[\frac{v}{2}]$.
In addition, from Observation \ref{obs}, it follows that $2x_2\leq
12$. Thus
$$\begin{array}{ll}
24-2v&=2x_4+4x_6+\cdots+(2k-2)x_{2k}\\
&\geq2(x_4+x_6+\cdots+x_{2k})\\
&=2(v-x_2)\geq 2v-12.
\end{array}$$
So $v\leq 9$ and $v=9$ if and only if $x_2=6$, $x_4=3$, and
$x_6=x_8=\ldots=x_{2k}=0$. So the degree sequence is
$<2,2,2,2,2,2,4,4,4>$ and graph should be of the following form:
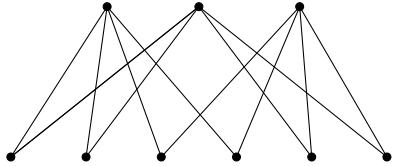
\begin{center}
\begin{tikzpicture}[line cap=round,line join=round,>=triangle 45,x=1.0cm,y=1.0cm]
\clip(-1.,-0.5) rectangle (6.,2.5);
\draw (1.28,2.)-- (0.,0.);
\draw  (1.28,2.)-- (1.,0.);
\draw  (1.28,2.)-- (2.,0.);
\draw  (1.28,2.)-- (3.,0.);
\draw  (2.5,2.)-- (0.,0.);
\draw  (2.5,2.)-- (1.,0.);
\draw  (2.5,2.)-- (4.,0.);
\draw  (2.5,2.)-- (5.,0.);
\draw  (3.84,2.)-- (5.,0.);
\draw (3.84,2.)-- (4.,0.);
\draw (3.84,2.)-- (3.,0.);
\draw  (3.84,2.)-- (2.,0.);
\draw  (2.5,2.)-- (0.,0.);
\begin{scriptsize}
\draw [fill=black] (0.,0.) circle (1.5pt);
\draw [fill=black] (1.,0.) circle (1.5pt);
\draw [fill=black] (2.,0.) circle (1.5pt);
\draw [fill=black] (3.,0.) circle (1.5pt);
\draw [fill=black] (4.,0.) circle (1.5pt);
\draw [fill=black] (5.,0.) circle (1.5pt);
\draw [fill=black] (2.5,2.) circle (1.5pt);
\draw [fill=black] (1.28,2.) circle (1.5pt);
\draw [fill=black] (3.84,2.) circle (1.5pt);
\end{scriptsize}
\end{tikzpicture}
\end{center}
which can be decomposed to $4-$cycles in a unique way.
\end{proof}

%\end{center}
\section{Generating $4-{\rm CS}(9)$'s from each other}

To generate two $4-{\rm CS}(9)$'s from each other by using  $4-$cycle bitrades of volume $2$ and $3$, first of all one needs
to be familiar with $4-{\rm CS}(9)$'s.
In \cite{MR1301216}, the authors characterized all non-isomorphic $4-{\rm CS}(9)$'s as follows.

\begin{figure}[!h]
    \centering
    \begin{tikzpicture}[line cap=round,line join=round,>=triangle 45,x=1.0cm,y=1.0cm]
    \clip(-0.7,-0.7) rectangle (3.7,3.7);
    \draw (1.,3.)-- (1.,0.);
    \draw (0.,2.)-- (1.,0.);
    \draw (0.,1.)-- (2.,0.);
    \draw (2.,0.)-- (1.,0.);
    \draw (1.,0.)-- (3.,1.);
    \draw (3.,1.)-- (3.,2.);
    \draw (2.,0.)-- (2.,3.);
    \begin{scriptsize}
    \draw [fill=black] (1.,3.) circle (2.0pt);
    \draw[color=black] (1.03,3.44) node {$S_1$};
    \draw [fill=black] (2.,3.) circle (2.0pt);
    \draw[color=black] (2.07,3.42) node {$S_2$};
    \draw [fill=black] (0.,2.) circle (2.0pt);
    \draw[color=black] (-0.43,2.12) node {$S_8$};
    \draw [fill=black] (0.,1.) circle (2.0pt);
    \draw[color=black] (-0.41,1.06) node {$S_7$};
    \draw [fill=black] (1.,0.) circle (2.0pt);
    \draw[color=black] (1.01,-0.26) node {$S_6$};
    \draw [fill=black] (2.,0.) circle (2.0pt);
    \draw[color=black] (2.03,-0.24) node {$S_5$};
    \draw [fill=black] (3.,1.) circle (2.0pt);
    \draw[color=black] (3.31,1.08) node {$S_4$};
    \draw [fill=black] (3.,2.) circle (2.0pt);
    \draw[color=black] (3.39,2.1) node {$S_3$};
    \end{scriptsize}
    \end{tikzpicture}
    \caption{A spanning tree.}\label{secret}
\end{figure}
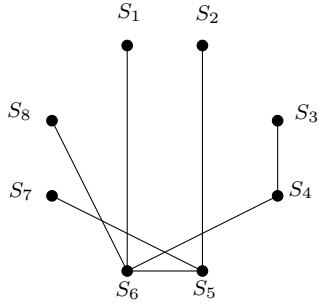

\begin{table}[h!]\begin{center}
 \begin{tabular}{c|ccccccccc}
$S_1$&1234&1356&1527&1829&2476&3648&3759&4589&6879\\
$S_2$&1234&1356&1527&1829&2476&3687&3849&4596&5798\\
$S_3$&1234&1356&1527&1829&2486&3647&3859&4579&6789\\
$S_4$&1234&1356&1527&1829&2486&3647&3879&4589&5769\\
$S_5$&1234&1356&1527&1829&2486&3678&3749&4596&5798\\
$S_6$&1234&1356&1527&1829&2486&3678&3759&4589&4697\\
$S_7$&1234&1356&1527&1829&2486&3698&3749&4576&5879\\
$S_8$&1234&1356&1527&1849&2458&2689&3678&3759&4697\\
\end{tabular}
\caption{All non-isomorphic  $4-{\rm CS}(9)$'s,~\cite{MR1301216}.}
\label{table4CS9}\end{center}
\end{table}
\begin{obs}\label{ptree}
    There exists a spanning tree corresponded to all $8$ non-isomorphic  $4-{\rm
    CS}(9)$s shown
    in  Table~{\rm\ref{table4CS9}} in which each vertex is corresponded to one of the given $4-{\rm
    CS}(9)$s, and each edge is corresponded to a  $4-$cycle bitrade of
    volume $3$ between two of its vertices.
\end{obs}
\begin{proof}
    A spanning tree is presented in Figure \ref{secret}. The vertices of tree are all eight
    non-isomorphic $4-{\rm CS}(9)$  mentioned in Table~\ref{table4CS9} and the edges
    are  $4-$cycle bitrades
    isomorphic to either $\mathcal{T}^{'}$ or to $\mathcal{T}^{''}$ in the proof of Theorem~\ref{Thm:Tv3f7}.
    The edges $e_{16},e_{34},e_{46},e_{56}$ and $e_{68}$ are
    corresponded
    to $\mathcal{T}^{'}$ and the edges $e_{25}$ and $e_{27}$ are corresponded  to $\mathcal{T}^{''}$.
\end{proof}

Let $\sigma$ be an arbitrary permutation on $\{1,\dots,9\}$ and $A$ be a subset of $4-$cycles in $K_9$. Then $A^{\sigma}$ is the subset
of  $4-$cycles  in $K_9$ obtained by applying $\sigma$ on each element of $A$. For example $\{(1,2,3,4)\}^{(1\ 5\ 3)}=\{(5,2,1,4)\}$. Thus, every $4-{\rm CS}(9)$ can be shown as $S_i^{\sigma}$, for some $S_i$  in Table~\ref{table4CS9} and some permutation $\sigma$ on $\{1,\dots,9\}$.

\begin{lem}\label{lem0}
If a  $4-{\rm CS}(9)$  $B$ can be generated from $4-{\rm CS}(9)$  $A$ by applying a sequence of trades $T_1,T_2,\ldots,T_k$, then
one can obtain $B^{\sigma}$ from $A^{\sigma}$ by applying $T_1^{\sigma},T_2^{\sigma},\ldots,T_k^{\sigma}$.
\end{lem}
Therefore by Observation~\ref{ptree}, we have the following corollary.

\begin{cor}\label{lem1}
       For every $S_i$ and $S_j$ in Table~$\ref{table4CS9}$ and every permutation $\sigma$ on $\{1,\dots,9\}$,
       there exists a path between $S_i^{\sigma}$ and $S_j^{\sigma}$ such that each edge is a bitrade of volume $3$.
\end{cor}
In addition, if there is a path between  $A$ and $A^{\sigma}$ for some permutation $\sigma$, then for each permutation $\tau$, there is a path between $A^{\tau}$ and $A^{\sigma\tau}$. This leads to  the following result.
\begin{cor}\label{cor2}
Let $A$ be a $4-{\rm CS}(9)$. If $A^{\sigma}$ and $A^{\tau}$ can be generated  from $A$ by applying
sequences of $4-$cycle trades of volume $2$ and $3$, then so is $A^{\sigma\tau}$.
\end{cor}
Suppose that $S_i^{\sigma}$ can be generated    from $S_i$ by applying a
sequence of $4-$cycle trades of volume $2$ and $3$ for some $i\in \{1,\ldots,9\}$ and a permutation $\sigma$ on $\{1,\dots,9\}$. Since by  Observation~\ref{ptree}, there exists a path of trades between  $S_i$ and $S_j$ and by Corollary~\ref{lem1}, one can generate $S_j^{\sigma}$ from $S_j$, for every $j\in \{1,\ldots,9\}$, thus, we have the following result.
\begin{lem}\label{lem2}
If there is a path of $4-$cycle bitrades of volume $2$ and $3$ between $S_i$ and $S_i^{\sigma}$ for some $i\in \{1,\ldots,9\}$ and a permutation $\sigma$ on $\{1,\dots,9\}$ then $S_j$ and $S_j^{\sigma}$ can be generated from each other by applying a
    sequence of  $4-$cycle bitrades of volume $2$ and $3$, for every $j\in \{1,\ldots,9\}$.
\end{lem}
In the following proposition we show that one can generate  $S_i^{(a\ 9)}$  from $S_j$ in Table~$\ref{table4CS9}$, for $a,i,j \in \{1,\ldots,8\}$.

\begin{prop}\label{prop:trp}
 $S_i^{(a\ 9)}$ can be generated from $S_j$ in Table~$\ref{table4CS9}$, for $a,i,j \in \{1,\ldots,8\}$,
    by applying a sequence of  $4-$cycle bitrades of volume $2$ and $3$.
\end{prop}

\begin{proof}
We show that for each $a\in\{1,\ldots,8\}$, the statement is true for $i=j=8$.
Then by Lemma~\ref{lem2} and Observation~\ref{ptree},  the statement  would be true for every $i$ and $j \in \{1,\ldots,8\}$.

 In the sequel, to generate $S_8^{(a\ 9)}$ from $S_8$, we consider the cycles which are different. For this purpose, we use three type $4-$cycle trades: $\mathcal{T}^{'}=(T^{'}_1,T^{'}_2),\
 \mathcal{T}^{''}=(T^{''}_1,T^{''}_2)$ and the double-diamond $D$.

First for moving from $S_1$ to $S^{(7\ 9)}_1$, we have
        $$\begin{array}{cccccc}
       \ \ \ \ S_1: &1527&1829&2476&4589&6879\\
        +\ \mathcal{T}^{''}: &&&{\bf 2496}&{\bf 4587}&{\bf 6897}\\
        +\ \ \mathcal{D}: &{\bf 1529}&{\bf 1827}&&&\\
        \end{array}$$
Thus, by Lemma \ref{lem2}, $S_8^{(7\ 9)}$ can be generated from $S_8$. Let $\sigma=(1\ 3\ 6\ 8\ 2\ 5\ 7\ 9\ 4)$. Since $S_8=S_8^{\sigma}$, we can consider a path between $S_8^{\sigma}$ and $S_8^{\sigma(7\ 9)}$ and by Lemma \ref{lem0}, there is a path between $S_8=S_8^{\sigma\sigma^{-1}}$ and  $S_8^{(4\ 9)}=S_8^{\sigma(7\ 9)\sigma^{-1}}$. Repeating this method, one can find  paths from $S_8$ to $S_8^{(4\ 1)}$, $S_8^{(1\ 3)}$, $S_8^{(3\ 6)}$, $S_8^{(6\ 8)}$, $S_8^{(8\ 2)}$ and $S_8^{(2\ 5)}$.

Note that $S_8^{(a\ 9)}=S_8^{(a\ b)(b\ 9)(a\ b)}$. So by Corollary \ref{cor2}, if we have  paths from $S_8$ to $S_8^{(a\ b)}$ and $S_8^{(b\ 9)}$, then we can find a path between $S_8$ and $S_8^{(a\ 9)}$.

Until now, we find pathes from $S_8$ to $S_8^{(7\ 9)}$ and $S_8^{(4\ 9)}$. From the path to $S_8^{(4\ 1)}$, we get a path to $S_8^{(1\ 9)}$. Then by the path to $S_8^{(1\ 3)}$, we get a path to $S_8^{(3\ 9)}$ and so on. This completes the proof.
\vspace{-.65cm}\end{proof}

It is known that every permutation of a finite set can be written as a product of transpositions. Even more, every permutation
can be written as a product of transposition in which all of them have a common element.
\begin{lem}\label{lem:per}
 If $a_j \neq a$ for $j \in \{1,\ldots,k\}$, then
$$(a_1\ a_2\ \ldots\ a_k)=(a_1\ a)(a_k\ a)\ldots(a_2\ a)(a_1\ a).$$
\end{lem}
Using this lemma and Corollary \ref{cor2}, we can extend Proposition \ref{prop:trp} to each arbitrary permutation $\sigma$.
\begin{thm}\label{Th:premain}
    For every $S_i$ and $S_j$  in Table~$\ref{table4CS9}$ and every permutation $\sigma$ on $\{1,\dots,9\}$,
       there exists a path between $S_j$ and $S_i^{\sigma}$ such that each edge is a bitrade of volume $2$ or $3$.
\end{thm}

\begin{thm}\label{Th:main}
    Every two $4-{\rm CS}(9)$s can be generated from each other by applying a
    sequence of  $4-$cycle bitrades of volume $2$ and $3$.
\end{thm}
\begin{proof}
Let $A$ and $B$ be two arbitrary $4-{\rm CS}(9)$s. There exist $i,j\in \{1,...,8\}$ and two
   permutations $\sigma$ and $\gamma$ on $\{1,\dots,9\}$   such that
   $A=S_i^{\sigma}$ and $B=S_j^{\gamma}$ where
    $4-$cycle systems $S_i^{\sigma}$ and $S_j^{\gamma}$ are obtained by applying permutation $\sigma$ and $\gamma$ on  $S_i$ and $S_j$ in Table~\ref{table4CS9}, respectively.
    By Theorem~\ref{Th:premain}, there exists a path of $4-$cycle trades of volume $2$ and $3$ from $S_i^{\sigma}$ to $S_j$ and
    there is a a path of $4-$cycle trades of volume $2$ and $3$ from $S_j$ to $S_j^{\gamma}$. The combination of these two paths completes the proof.
\end{proof}

%--------------------------------------------------------------------------------------------------------------------------

Theorem~\ref{Th:main} failed to be true for arbitrary $n$. In fact there are some $4-{\rm CS}(n)$ for proper $n$ with no double-diamond or any structure of $T'_i$ or $T_i''$.
\begin{example}
    The following cycles generate two cyclic $4-{\rm CS}(25)$ and $4-{\rm CS}(49)$, respectively, which don't contain any double-diamond or any structure of $T'_i$ or $T_i''$.\\
    {\bfseries{$4-{\rm CS}(25):$}}
    $(0 , 3, 1, 12)\ (0 , 4, 10, 17)\ (0 , 1, 6, 15)$\\
    {\bfseries{$4-{\rm CS}(49):$}}
    $(0 ,23, 20, 29)\   (0 , 22, 17, 32)\   (0 , 4, 12, 37)$
    $(0 ,18, 46, 36)\   (0 , 30, 44, 11)\   (0 , 1, 8, 2)$
\end{example}

\section{Linear algebraic presentation }

Let $M$ be a pair inclusion matrix whose
rows are corresponded to the  edges of the complete graph $K_n$
and its columns are corresponded to  all possible $4-$cycles of the complete graph $K_n$. An entry  $M_{e,C}$ is 1, if the edge $e$ belong to the cycle $C$. It is 0 otherwise.

Since by every four vertices $a, b, c$ and $d$, we can construct
three  different $4-$cycles $(a,b,c,d)$, $(a,c,b,d)$ and
$(a,b,d,c)$, the matrix $M$ has exactly $3{n\choose 4}$ columns.
Thus, the matrix
$M$ is of size ${n\choose 2}\times3{n\choose 4} $.

Also, for each trade $(T_1,T_2)$, we  consider a ``frequency''
vector $X$ with $3{n\choose 4}$ components with 1 for each cycle
in $T_1$ and $-1$ for each cycle in $T_2$, other cycles are
corresponded  with a 0  component. It is easy to see that every
vector $X$ corresponding to a trade is a vector in the kernel of
$M$. 

\begin{thm}
    The pair inclusion matrix $M$ is a full rank matrix.
\end{thm}
\begin{proof}
    Let $R_i$ be the row of matrix $M$ corresponding to the edge $e_i$
    and let  $\sum_{i=1}^{n\choose2} \lambda_iR_i=\bar{0}$ for some
    $\lambda_i\in\mathbb{R}$, where $\bar{0}$ is a $1\times 3{n\choose 4}$ vector with all entry are equal to $0$. First we show that for each two adjacent
    edges $e_i,e_j$, we have $\lambda_i+\lambda_j=0$. Consider cycles
    in Figure \ref{two}.
    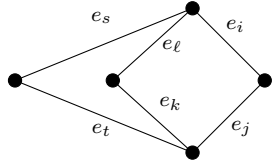
\begin{figure}[!h]
        %   \vspace{-2cm}
        \centering

        \begin{tikzpicture}[line cap=round,line join=round,>=triangle 45,x=1.0cm,y=1.0cm]
        \clip(0.,0.) rectangle (4.,3.);
        \draw (2.86,2.62)-- (3.82,1.66);
        \draw (3.82,1.66)-- (2.86,0.7);
        \draw (2.86,0.7)-- (1.8,1.66);
        \draw (1.8,1.66)-- (2.86,2.62);
        \draw (2.86,2.62)-- (0.5,1.66);
        \draw (0.5,1.66)-- (2.86,0.7);
        \begin{scriptsize}
        \draw [fill=black] (2.86,2.62) circle (2.5pt);
        \draw [fill=black] (3.82,1.66) circle (2.5pt);
        \draw[color=black] (3.43,2.38) node {$e_i$};
        \draw [fill=black] (2.86,0.7) circle (2.5pt);
        \draw[color=black] (3.51,1.02) node {$e_j$};
        \draw [fill=black] (1.8,1.66) circle (2.5pt);
        \draw[color=black] (2.57,1.36) node {$e_k$};
        \draw[color=black] (2.59,2.12) node {$e_{\ell}$};
        \draw [fill=black] (0.5,1.66) circle (2.5pt);
        \draw[color=black] (1.65,2.50) node {$e_s$};
        \draw[color=black] (1.65,1.0) node {$e_t$};
        \end{scriptsize}
        \end{tikzpicture}
        %\vspace{-3cm}
        \caption{Two $4-$cycles.}\label{two}
    \end{figure}

    From $\sum_{i=1}^{n\choose2} \lambda_iR_i=0$, it follows that
    $$
    \begin{array}{c}
    \lambda_i+\lambda_j+\lambda_k+\lambda_{\ell}=0,\\
    \lambda_i+\lambda_j+\lambda_s+\lambda_{t}=0,\\
    \lambda_s+\lambda_t+\lambda_k+\lambda_{\ell}=0.
    \end{array}
    $$
    Thus, $\lambda_i+\lambda_j=\lambda_k+\lambda_{\ell}=\lambda_s+\lambda_t$ and therefore $\lambda_i+\lambda_j=0$.

    Now, if three edges $e_i,e_j,e_r$ form a triangle, we have
    $$
    \begin{array}{c}
    \lambda_i+\lambda_j=0,\\
    \lambda_i+\lambda_r=0,\\
    \lambda_r+\lambda_j=0.
    \end{array}
    $$
    Therefore $\lambda_i=0$. That is, the columns of $M$ are linearly independent.
\end{proof}
\begin{cor}
    The nullity of $M$ is $3{n\choose 4}-{n\choose 2}$.
\end{cor}
{\bf Open problem}
The set of vectors corresponding with the
double-diamonds
    is a generating set for the kernel of $M$. In other words, for
    each $n$, there exists a set of $3{n\choose 4}-{n\choose 2}$
    linearly independent vectors in the kernel of $M$, where each of
    them is corresponded to a double-diamond.

\section*{Statements \& Declarations}
The authors declare that no funds, grants, or other support were received during the preparation of this manuscript.

The authors have no relevant financial or non-financial interests to disclose.

All authors contributed to the design and implementation of the research, to the analysis of the 
results and to the writing of the manuscript. All authors read and approved the final manuscript.

\begin{thebibliography}{20}
  \bibitem{MR1821945}
  {P. Adams,  E. J. Billington, D. E. Bryant and A. Khodkar},
  {The $\mu-$way intersection problem for $m-$cycle systems}, \textit{Discrete Math.}, {\bf 231} (2001), 27--56.
  % % % % % % % % % % % % % % % % % % % % % % % % % % % % % % % % % % % % % % % %
  \bibitem{MR2195316}
  P. Adams, D. Bryant, M. Grannell and T. Griggs, Diagonally
  switchable $4-$cycle systems, \textit{ Australas. J. Combin.}, {\bf
  34} (2006), 145--152.
  % % % % % % % % % % % % % % % % % % % % % % % % % % % % % % % % % % % % % % % %
  \bibitem{MR2883525}
  {M. Aryapoor and E. S. Mahmoodian}, {On uniformly generating {L}atin squares}, \textit{Bull. Inst. Combin. Appl.}, {\bf 62} (2011), 48--58.
  % % % % % % % % % % % % % % % % % % % % % % % % % % % % % % % % % % % % % % % % %
  \bibitem{MR1392993}
  C. J. Colbourn and J. H. Dinitz, \textit{The CRC handbook of
      combinatorial designs}, CRC Press, 2010.
      % % % % % % % % % % % % % % % % % % % % % % % % % % % % % % % %
  \bibitem{MR1301216}
I. Dejter,  P. Rivera-Vega and A. Rosa, {Invariants for $2-$factorizations and cycle systems},
\textit{J. Combin. Math. Combin. Comput.}, {\bf 16} (1994), 129--152.


\bibitem{MR1410617}
        M. T. Jacobson and P. Matthews, Generating uniformly
        distributed random Latin squares, \textit{J. Combin. Des.}, {\bf
        4(6)} (1996),  405--437.
   % % % % % % % % % % % % % % % % % % % % % % % % % % % % % % % % %
      \bibitem{MR2976352}
          A. A. Khanban, M. Mahdian and E. S.  Mahmoodian, A linear algebraic
          approach to orthogonal arrays and Latin squares, \textit{Ars
          Combin.}, {\bf105} (2012), 15--22.

  \bibitem{MR1954532}
   {C. C. Lindner}, { A partial {$4-$}cycle system of order {$n$}
   can be embedded in a {$4-$}cycle system of order at most {$2n+15$}},
  \textit{Bull. Inst. Combin. Appl.}, {\bf 37} (2003), 88--93.

  \bibitem{MR1367739}
    {D. B. West}, {Introduction to graph theory},
    {Prentice Hall, Inc., Upper Saddle River, NJ},
  {2001}.
%%%%%%%%%%%%%%%%%%%%%%%%%%%%%%%%%%5


\end{thebibliography}
 \end{document}